\documentclass[a4paper,reqno]{amsart}
\usepackage{amssymb}

\newenvironment{Proof}{{\sc Proof.}\ }{~\rule{1ex}{1ex}\vspace{0.5truecm}}
\newtheorem{Lemma}{Lemma}[section]
\newtheorem{Th}[Lemma]{Theorem}
\newtheorem{Prop}[Lemma]{Proposition}
\newtheorem{Rem}[Lemma]{Remark}
\newtheorem{Cor}[Lemma]{Corollary}

\begin{document}
\title[]{Iterated power intersections of ideals in rings of iterated differential polynomials II}
\author{Pavel P\v r\'\i hoda}
\address{Charles University, Faculty of Mathematics and Physics,
Department of Algebra, Sokolovsk\' a 83, 18675 Prague 8, Czech Republic}
\email{prihoda@karlin.mff.cuni.cz}
\thanks{This research was supported by Czech Science Foundation grant GA \v CR 20-13778S}
\keywords{}
\subjclass[2010]{16D25,16D40,16S30}
\begin{abstract}
  This paper is a continuation of \cite{PP} where iterated intersections of powers of ideals were studied in 
	rings of iterated differential polynomials. We present a method which can be used to show that for every proper 
	ideal $I$ of a suitable ring of iterated differential polynomials almost all iterated intersections of powers 
	of $I$ have to be zero.
	The setback of the method is that it works only if the derivations used in the construction of the ring of 
	iterated differential polynomials satisfy additional assumptions. On the other hand, it can be applied in cases 
	which were not covered by \cite{PP}, for example for a universal enveloping algebra of a completely 
	solvable Lie algebra over a field of positive characteristic.
	
\end{abstract}
\maketitle
\section{Introduction}

Let $I$ be an ideal of a ring $R$ and let $I(1):= \cap_{n \in \mathbb{N}} I^n$ denote the intersection of all its powers.
Let us define inductively $I(m):=(I(m-1))(1)$ for every $2\leq m \in \mathbb{N}$. The ideal $I(m)$ is called the {\em $m$-iterated 
intersection of powers of $I$}. This short note gives some examples of rings 
having the following property: For any proper ideal $I$ of $R$ there exists $m \in \mathbb{N}$ such that $I(m) = 0$.

This property was studied for rings of iterated differential polynomials in \cite{PP}. Recall that if $R$ is a ring 
and $\delta\colon R \to R$ is a {\em derivation} on $R$, i.e., an additive map satisfying $\delta(rs) = \delta(r)s+r\delta(s)$
for every $r,s \in R$, the {\em ring of differential polynomials} $R[x,\delta]$ is a free left $R$-module with basis $1,x,x^2,\dots$
with a multiplication given by the rule $xr - rx = \delta(r), r \in R$ (see for example \cite[Chapter~2]{GW} for details of this construction).
Iterating this construction $n$-times we get a ring $S:=R[x_1,\delta_1,\dots,x_n,\delta_n]$ of {\em iterated differential polynomials over $R$}.
For example, if $R_i := R[x_1,\dots,x_i]$ is the subring of $S$ generated by $R$ and $x_1,\dots,x_i$, then $\delta_{i+1}$ is a derivation 
on $R_i$ and $R_{i+1} := R[x_1,\dots,x_{i+1}]$ is isomorphic to $R_i[x,\delta_{i+1}]$. 

By \cite[Theorem~3.4]{PP}, if $R$ is a commutative noetherian domain which is a $\mathbb{Q}$-algebra and $I$ is a proper 
ideal of $R[x_1,\delta_1,\dots,x_n,\delta_n]$ a ring of iterated differential polynomials over $R$, then $I(n+1) = 0$.

A natural source of examples of rings of differential polynomials are universal enveloping algebras of 
solvable Lie algebras of finite dimension. By \cite[Proposition~4.1]{PP}, if $L$ is a solvable Lie algebra  
of finite dimension over a field of characteristic $0$, then $I(2) = 0$ for every proper ideal $I$ of the universal enveloping 
algebra of $L$.

The approach used in \cite{PP} depends on the fact that each prime ideal in the considered ring of iterated differential polynomials
is completely prime (cf. \cite[Theorem~10.23]{GW} which is a key ingredient in the proof of \cite[Theorem~3.4]{PP}). 
In particular, this approach does not work for algebras over fields of positive characteristic.

We offer an elementary method which can be applied to a differential polynomial ring 
$S = R[x_1,\delta_1,x_2,\delta_2,\dots,x_n,\delta_n]$ even in positive characteristic provided the derivations satisfy additional properties.
For example, assume that
\begin{enumerate}
\item[(a)] $\delta_i(R) \subseteq R$ for any $1 \leq i \leq n$,
\item[(b)] for every $1 \leq j < i \leq n$ there are $u,v \in R[x_1,\dots,x_{j-1}]$ such that $\delta_i(x_j) = ux_j + v$,
\item[(c)] $R$ is right noetherian and there are no ideals of $R$ invariant under $\delta_1,\dots,\delta_n$ with exception of $0$ and $R$.
\end{enumerate}
In Theorem \ref{t2} we prove that in this case every proper ideal $I$ of $S$ satisfies $I(n) = 0$. A natural example of a ring of differential 
polynomials satisfying $(a),(b)$ and $(c)$ is a universal enveloping algebra of a completely solvable Lie algebra over a field. 
 
Our main motivation to study iterated power intersections of ideals stems from their relation with the absence of interesting non-finitely 
generated projective modules.

\begin{Prop} \label{p1} (\cite[Proposition~3.7]{PP}) Let $S$ be a ring such that $S/J(S)$ is right noetherian. 
Assume that for every proper ideal $I$ of $S$ there exists $m \in \mathbb{N}$
such that $I(m) = 0$. Then every projective right $S$-module which is not finitely generated is free. 
\end{Prop} 
\begin{Proof}
As explained in the proof of \cite[Proposition~3.7]{PP}, the assumption on ideals of $S$ implies that every countably but not finitely generated  
projective right $S$-module $P$ is uniformly $\aleph_0$-big in the sense of Bass \cite{B}, that is, if $P/PI$ is finitely generated for some ideal $I$ of $S$,
then $I = S$. The result then follows from \cite[Theorem~3.1]{B}.
\end{Proof}

\section{Reducing the number of variables via localization and factorization} \label{s2}

Any ring of iterated differential polynomials $S = R[x_1,\delta_1,\dots,x_n,\delta_n]$ in $n$ variables over $R$ can be considered 
as a ring of iterated differential polynomials in $n-1$ variables over $R[x_1,\delta_1]$, i.e., $S = (R[x_1,\delta_1])[x_2,\delta_2,\dots,x_n,\delta_n]$.
In the proofs of Theorem~\ref{t1} and Theorem~\ref{t2} we consider right localizations of $S$ with respect to right denominator sets 
of $R$ and $R[x_1,\delta_1]$ and factors of $S$ modulo ideals of $S$ which are generated by ideals of $R[x_1,\delta_1]$ invariant under derivations $\delta_2,\dots,\delta_n$.
Although these constructions are quite standard, in this section we recall them in detail in order to explain the notation which is 
used in the proofs of the results in Section \ref{s3}. 

Let $R$ be a ring, $\delta\colon R \to R$ a derivation on $R$ and let $X \subseteq R$ be a right denominator set of $R$ (for details 
on right quotient rings see for example \cite[Chapter~2]{MR}).
By \cite[Exercise 10R]{GW}, $X$ is a right denominator set of $S = R[x,\delta]$ and the right localization $SX^{-1}$ of $S$ with respect to $X$
has a structure of a differential polynomial ring. To be more precise, let $\theta \colon R \to RX^{-1}$ be a right localization 
of $R$ with respect to $X$. Then there exists a unique derivation $\partial \colon RX^{-1} \to RX^{-1}$ on $RX^{-1}$ 
satisfying $ \partial \theta = \theta \delta$ and the canonical map 
$\varphi \colon R[x,\delta] \to RX^{-1}[x,\partial]$ which maps $\sum_{i = 0}^{d} r_ix^i$ to $\sum_{i = 0}^d \theta(r_i)x^i$ is 
a right localization of $R[x,\delta]$ with respect to $X$.

Of course, this construction can be iterated: If $S = R[x_1,\delta_1,\dots,x_n,\delta_n]$ is a ring of iterated differential polynomials
and $X\subseteq R$ is a right denominator set of $R$, then $X$ is a right denominator set in $S$ and the ring $SX^{-1}$ has a 
structure of a ring of iterated differential polynomials. Let $\theta \colon R \to RX^{-1}$ be as above and suppose we have defined 
a ring $RX^{-1}[x_1,\partial_1,\dots,x_i,\partial_i]$ (i.e., we have defined derivations $\partial_1,\dots,\partial_i$) such that the 
canonical map $$\varphi_i \colon R[x_1,\delta_1,\dots,x_i,\delta_i] \to RX^{-1}[x_1,\partial_1,\dots,x_i,\partial_i]$$ 
which sends a term $rx_1^{e_1}\cdots x_i^{e_i}$ to $\theta(r)x_1^{e_1}\cdots x_i^{e_i}$ is a right localization 
of $R[x_1,\delta_1,\dots,x_i,\delta_i]$ with respect to $X$. By the single variable case discussed above, $X$ is a right denominator set in 
$R[x_1,\delta_1,\dots,x_{i+1},\delta_{i+1}]$ and there exists a unique derivation 
$$\partial_{i+1}\colon RX^{-1}[x_1,\partial_1,\dots,x_{i},\partial_{i}] \to RX^{-1}[x_1,\partial_1,\dots,x_{i},\partial_{i}]$$ 
on $RX^{-1}[x_1,\partial_1,\dots,x_{i},\partial_{i}]$ such that 
$\varphi_{i}\delta_{i+1} = \partial_{i+1} \varphi_i$. Then the canonical map $$\varphi_{i+1} \colon R[x_1,\partial_1,\dots,x_{i+1},\partial_{i+1}] 
\to RX^{-1}[x_1,\partial_1,\dots,x_{i+1},\partial_{i+1}]$$ which sends a term $rx_1^{e_1}\cdots x_{i+1}^{e_{i+1}}$ to $\theta(r)x_1^{e_1}\cdots x_{i+1}^{e_{i+1}}$
is a right localization of $R[x_1,\delta_1,\dots,x_{i+1},\delta_{i+1}]$ with respect to $X$.
 
Repeating this process $n$-times, we define a ring $RX^{-1}[x_1,\partial_1,\dots,x_n,\partial_n]$ such that the canonical map 
$$\varphi \colon S \to RX^{-1}[x_1,\partial_1,\dots,x_n,\partial_n] $$ which sends a term $rx_1^{e_1}\cdots x_n^{e_n}$ to 
$\theta(r)x_1^{e_1}\cdots x_n^{e_n}$ is a right localization of $S$ with respect to $X$. 

Now consider a ring of iterated differential polynomials $S = R[x_1,\delta_1,\dots,x_n,\delta_n]$ and an ideal 
$P \subseteq R$ of $R$ such that $\delta_i(P) \subseteq P$ for every $1 \leq i \leq n$. Then $PS = P[x_1,x_2,\dots,x_n]$ is the ideal of $S$ generated by $P$.
The factor ring $\overline{S} = S/PS$ can be seen as a ring of iterated differential polynomials over $R/P$:
Let $\pi \colon R \to R/P$ be the canonical projection. Since $P$ is $\delta_1$-invariant, there exists a derivation $\overline{\delta_1} \colon R/P \to R/P$
such that $\pi \delta_1 = \overline{\delta_1}\pi$. The canonical map $\pi_1 \colon R[x_1,\delta_1] \to R/P[x_1,\overline{\delta_1}]$ is a homomorphism 
of rings with kernel $P[x_1]$.

Assume we have defined $\overline{\delta_1},\dots,\overline{\delta_{i}}$ such that the canonical map
$$\pi_i \colon R[x_1,\delta_1,\dots,x_i,\delta_i] \to R/P[x_1,\overline{\delta_1},\dots,x_i,\overline{\delta_i}]$$
which sends $rx_1^{e_1}\cdots x_i^{e_i}$ to $\pi(r)x_1^{e_1}\cdots x_i^{e_i}$ is an onto ring homomorphism with kernel $P[x_1,\dots,x_i]$.
Then $P[x_1,\dots,x_i]$ is stable under $\delta_{i+1}$. Then there exists a derivation $\overline{\delta_{i+1}} \colon R/P[x_1,\overline{\delta}_1,\dots,x_{i},\overline{\delta}_{i}]
\to R/P[x_1,\overline{\delta}_1,\dots,x_{i},\overline{\delta}_{i}]$ such that $\overline{\delta_{i+1}}\pi_i =\pi_i \delta_{i+1}$. Then the canonical map
$$\pi_{i+1} \colon R[x_1,\delta_1,\dots,x_{i+1},\delta_{i+1}] \to R/P[x_1,\overline{\delta_1},\dots,x_{i+1},\overline{\delta_{i+1}}]$$
is a ring homomorphism with kernel $P[x_1,\dots,x_{i+1}]$.

After $n$ repetitions of this process we define on $S/PS$ a structure of a ring of iterated differential polynomials.

\section{The results} \label{s3}

Let us single out several easy observations we are going to use repeatedly in the proofs.

\begin{Lemma} \label{l1}
Let $\varphi \colon R \to S$ be a homomorphism of rings, $I$ an ideal of $R$, $J$ an ideal of $S$.
If $\varphi(I) \subseteq J$ and $J(n) = 0$ for some $n \in \mathbb{N}$, then $\varphi(I(n)) = 0$.
\end{Lemma}

\begin{Proof}
If $\varphi(I) \subseteq J$, then it is easy to check by induction that $\varphi(I(k)) \subseteq J(k)$ for every $k \in \mathbb{N}$.
Then the result is immediate.
\end{Proof}

\begin{Lemma} \label{l2}
Let $R$ be a ring and let $I \subseteq R$ be its ideal satisfying $I(1) = 0$. If there exists $n \in \mathbb{N}$ such that 
every proper ideal $J$ of $R/I$ satisfies $J(n) = 0$, then for every ideal $L$ of $R$ is either $L + I = R$ or $L(n+1) = 0$.
\end{Lemma}

\begin{Proof}
Apply Lemma \ref{l1} to the canonical projection $\pi \colon R \to R/I$. If $L$ is an ideal of $R$ and $\pi(L) \neq R/I$, then $\pi(L(n)) = 0$, i.e., 
$L(n) \subseteq I$. Therefore $L(n+1) \subseteq I(1) = 0$.
\end{Proof}

\begin{Lemma} \label{l5}
Let $R$ be a ring and let $M_1,\dots,M_k,I$ be ideals of $R$ such that $(I + \cap_{j = 1}^k M_j)/I$ is a nil ideal of $R/I$.
If $I \neq R$, then there exists $1 \leq j \leq k$ such that $M_j + I \neq  R$. 
\end{Lemma}

\begin{Proof}
Assume $M_j + I = R$ for every $1 \leq j \leq k$. Then 
$$R = (I+M_1)\cdots (I+M_k) = I + M_1 \cdots M_k\,.$$
In particular, $1 = i + m$ for some $i \in I$ and $m \in \cap_{j = 1}^k M_j$. By our assumption, $m+I$ is nilpotent in $R/I$ and hence $I = R$.
\end{Proof}

If $S$ is a ring and $\Delta$ is a set of derivations on $S$, then a set $X \subseteq S$ is called {\em $\Delta$-invariant} if 
$\partial(X) \subseteq X$ for every $\partial \in \Delta$. If $\Sigma$ is a set of ideals of $S$, {\em maximal elements of $\Sigma$}
denotes the set of maximal elements of a partially ordered set $(\Sigma\setminus \{S\}, \subseteq)$.

\begin{Lemma} \label{l3}
Let $R$ be a ring, $\delta \colon R \to R$ a derivation on $R$, and $S := R[x,\delta]$ the corresponding ring of differential polynomials.
Let $\Delta$ be a set of derivations on $S$ satisfying the following conditions:
\begin{enumerate}
  \item[(a)] $R$ is $\Delta$-invariant,
  \item[(b)] $\partial(x)$ is a polynomial of degree at most $1$ for every $\partial \in \Delta$,
  \item[(c)] if $J \subseteq R$ is a $\Delta$-invariant ideal of $R$ such that $\delta(J) \subseteq J$, then either $J = 0$ or $J = R$.
\end{enumerate}
Let $I$ be a nonzero proper $\Delta$-invariant ideal of $S$. Then 
\begin{enumerate}
  \item[(i)] There exists unique $f \in Z(S)$ monic such that $I = fS = Sf$.
	\item[(ii)] If $R$ is right (left) artinian, then $S/I$ is right (left) artinian.
	\item[(iii)] Let $\Sigma$ be the set of all $\Delta$-invariant ideals of $S$ containing $I$. Then the set of maximal elements of $\Sigma$ 
							 is finite and if $M_1,\dots,M_k$ is a list of these elements, then $I = M_1^{e_1}\cdots M_k^{e_k}$ for some $e_1,\dots,e_k \in \mathbb{N}$.   
\end{enumerate}
\end{Lemma}

\begin{Proof}
Let $I$ be a nonzero proper $\Delta$-invariant ideal of $S$ and let $d := {\rm min}\{{\rm deg}(f) \mid 0 \neq f \in I\}$.  
Since $I \cap R$ is $\Delta$-invariant by (a) and $\delta(I \cap R) \subseteq I \cap R$, (c) implies $d > 0$.

(i) (cf. \cite[Proposition~2.1]{GW})
Consider the set $J := \{r \in R \mid \exists r_0,\dots,r_{d-1} \in R$ such that $rx^d + \sum_{i = 0}^{d-1} r_ix^i \in I\}$. Then $J$ is easily checked to be 
an ideal of $R$. If $f = \sum_{i = 0}^d r_ix^i \in I$ and $\partial \in \Delta$, then, regarding (b), $\partial(f) = \partial(r_d) x^d + r_d \partial(x^d) + g$, for some $g \in S$
of degree at most $d-1$. So the coefficient by $x^d$ in $\partial(f)$ is of the form $\partial(r_d) + r_ds$ for some $s \in R$. 
Since $\partial(f) \in I$ and $r_d \in J$, $\partial(r_d) \in J$. That is, $J$ is invariant under $\partial$. 
Also note that $xf - fx = \sum_{i = 0}^d \delta(r_i)x^i \in I$, thus $\delta(r_d) \in J$. By (c), $J = R$.

Therefore $I$ contains a monic polynomial $f$ of degree $d$. Moreover, it is easy to see that $I$ is generated by $f$ as a left ideal. 
Observe that polynomial $xf - fx$ has degree less than $d$ and is contained in $I$. It follows that $xf - fx = 0$, i.e, $f$ commutes with $x$. 
Similar argument shows that $f$ commutes with every $r \in R$, hence $f \in Z(S)$. 
Since $I$ contains no nonzero polynomials of degree less than $d$, $f$ is the only monic polynomial of degree $d$ contained in $I$.

(ii) Assume that $R$ is left (right) artinian.
Note that $x^i + I, i = 0,\dots,d-1$ generate $S/I$ as a left (right) $R$-module, so $S/I$ is left (right) artinian.

(iii) Let $\Sigma_m$ be the set of all maximal elements of $\Sigma$. If $M \in \Sigma_m$, then, by (i), there exists a unique monic $f_M \in M$ such 
that $M = f_MS$. Note that if $M \varsubsetneq N \in \Sigma$, then ${\rm deg}(f_M) > {\rm deg}(f_N)$. In particular, $\Sigma$ satisfies a.c.c.
We claim that $f:=f_I$ is a product of elements from the set $\{f_M \mid M \in \Sigma_m\}$.

Consider an equality $f = f_{M_1}\cdots f_{M_k} g$ with $M_1,\dots,M_k \in \Sigma_m$ and $g \in S$ of positive degree. Observe that 
$f_{M_1},\dots,f_{M_k},g$ are monic and $f,f_{M_1},\dots,f_{M_k} \in Z(S)$. Since monic polynomials are regular 
elements of $S$, $g \in Z(S)$. Let $\partial \in \Delta.$ Since $I,M_1,M_2,\dots,M_k$ are $\Delta$-invariant, there are $u,v_1,v_2,\dots,v_k \in S$ such that 
$\partial(f) = fu$, $\partial(f_{M_i}) = f_{M_i}v_i$ for $i = 1,\dots,k$. Apply $\partial$ on $f = f_{M_1}\cdots f_{M_k}g$:
$$fu = \partial(f) = \sum_{i = 1}^k f_{M_1}\cdots f_{M_{i-1}}(f_{M_i}v_i)f_{M_{i+1}}\cdots f_{M_k}g+  f_{M_1}\cdots f_{M_k} \partial(g)$$
$$fu = fh + f_{M_1}\cdots f_{M_k} \partial(g),$$
where $h = \sum_{i = 1}^k v_i$. Cancellation of the regular element $f_{M_1}\cdots f_{M_k}$ from the equality above gives
$\partial(g) = g(u - h)$. In particular, the ideal $gS$ is invariant under $\partial$ and hence $gS \in \Sigma$.
Then there exists $M_{k+1} \in \Sigma_m$ containing $g$ and $g = f_{M_{k+1}}g'$ for some $g' \in S$ with degree less than ${\rm deg}(g)$.
After finitely many repetitions of this process we get $f$ as a product of elements from the set $\{f_M \mid M \in \Sigma_m\}$.

Let us change the notation and write $f = f_{M_1}^{e_1}\cdots f_{M_k}^{e_k}$, where $M_1,\dots,M_k$ are pairwise different elements of $\Sigma_m$ and
$e_1,\dots,e_k$ are positive integers. Then $fS = M_1^{e_1}\cdots M_k^{e_k}$. Let us verify that $\{M_1,\dots,M_k\} = \Sigma_m$.
Assume that there exists $M \in \Sigma$ which is not contained in any element from $\{M_1,\dots,M_k\}$. Since $M+M_i$ is $\Delta$-invariant for every $1 \leq i \leq k$, 
the maximality of $M_i$ implies $M + M_i = S$. On the other hand, if $M' := \cap_{i = 1}^k M_i$, then $M'^e \subseteq fS$ for any $e \geq \sum_{i = 1}^k e_i$.
This implies that $M' + M /M$ is a nil ideal of $S/M$ and therefore $S = M$ by Lemma~\ref{l5}.
This  shows that $\Sigma_m = \{M_1,\dots,M_k\}\,.$
\end{Proof}

\begin{Th} \label{t1}
Let $R$ be a noetherian prime ring such that the abelian group $(R,+)$ is torsion free. Consider a ring of iterated differential polynomials 
$S = R[x_1,\delta_1,x_2,\delta_2,\dots,x_n,\delta_n]$ which satisfies: 
\begin{enumerate}
\item[(a)] $\delta_i(R) \subseteq R$ for every $i = 1,\dots,n$.
\item[(b)] $\delta_i(x_j) \in R[x_1,\dots,x_j]$ for every $1 \leq j < i \leq n$. 
\item[(c)] If $J$ is an ideal of $R$ such that $\delta_i(J) \subseteq J$ for every $i = 1,\dots,n$, then $J =0$ or $J = R$.
\end{enumerate}
Then $I(n) = 0$ for every proper ideal $I$ of $S$.
\end{Th}

\begin{Proof}
We claim that it is enough to consider the case when $R$ is a simple artinian ring of characteristic zero (and, in particular, a $\mathbb{Q}$-algebra).
Indeed, let $X$ be the set of all regular elements of $R$. Then, by \cite[Theorem 2.3.6]{MR}, $X$ is a right denominator set in $R$. Following the notation from 
Section \ref{s2}, $X$ is also a right denominator set in $S$  and the map 
$\varphi\colon S \to RX^{-1}[x_1,\partial_1,\dots,x_n,\partial_n] = SX^{-1}$ given by 
$$\varphi(r) = r.1^{-1},r \in R, \varphi(x_i) = x_i, i = 1,\dots,n$$
is a right localization of $S$ with respect to $X$. By Goldie's theorem \cite[Theorem 6.18]{GW}, $RX^{-1}$ is simple artinian and by 
Hilbert's basis theorem \cite[Theorem~2.6]{GW}, $SX^{-1}$ is left and right noetherian. Therefore it is possible to apply \cite[Theorem~10.18(a)]{GW}
to see that $\varphi(K)X^{-1}$ is an ideal of $SX^{-1}$ for every ideal $K$ of $S$.
Also note that $(R,+)$ is a torsion free abelian group, so 
$SX^{-1}$ is a $\mathbb{Q}$-algebra.
Assume that every proper ideal $U$ of $SX^{-1}$ satisfies $U(n) = 0$. Consider a proper ideal $I \subseteq S$. By (c), $I \cap R = 0$, in particular 
$I \cap X = \emptyset$. Then $\varphi(I)X^{-1}$ is a proper ideal of $SX^{-1}$ and, by Lemma \ref{l1}, $I(n) = 0$.  

In the rest of the proof we assume that $R$ is a simple artinian $\mathbb{Q}$-algebra. We proceed by induction on $n$.
Let $n = 1$ and consider a nonzero proper ideal $I \subseteq S = R[x_1,\delta_1]$. By Lemma \ref{l3}(i) applied for $\Delta = \emptyset$, $I$ is generated by
a monic polynomial from $Z(S)$. Then $I(1) = 0$.  

Now assume that the statement is true for rings with $n-1$ variables and consider a proper ideal $I \subseteq S = R[x_1,\delta_1,\dots,x_n,\delta_n]$.
Let $J: = I \cap R[x_1]$. We distinguish two cases:

1) $J = 0$. By \cite[Theorem 1.2.9]{MR}, $R[x_1] = R[x_1,\delta_1]$ is a noetherian prime ring, hence if $X$ the set of all its regular elements $R[x_1]X^{-1}$ is a simple artinian ring.
As in the first step of the proof, $X$ is a  right denominator set of $S$ and the homomorphism $\varphi \colon S \to (R[x_1]X^{-1})[x_2,\partial_2,\dots,x_n,\partial_n]$
given by $$\varphi(p) = p.1^{-1},p \in R[x_1], \varphi(x_i) = x_i, i = 2,\dots,n$$
is a right localization of $S$ with respect to $X$. Since $I \cap X = \emptyset$, $\varphi(I)X^{-1}$ is a proper ideal of $SX^{-1}$ and $\varphi(I)X^{-1}(n-1) = 0$ by 
the induction hypothesis. Then $I(n-1) = 0$ by Lemma \ref{l1}.

2) $0 \neq J$. Observe that $J$ is a proper ideal of $R[x_1]$. By (a) and (b), the restriction of $\delta_i$ to $R[x_1]$ is a ($\mathbb{Q}$-linear) derivation
of $R[x_1]$ and $\delta_i(J) \subseteq J$ for every $2 \leq i \leq n$. Let $P_1,\dots,P_k$ be a list of all primes of $R[x_1]$ which are minimal over $J$.
By \cite[Proposition~14.2.3]{MR}, all these primes are  stable under $\delta_2,\dots,\delta_n$. As explained in Section \ref{s2}, for each $i = 1,\dots,k$ the ring 
$S/P_iS$ has a structure of a ring of iterated differential polynomials,  $S/P_iS \simeq (R[x_1]/P_i)[x_2,\overline{\delta_2},\dots,x_n,\overline{\delta_n}]$.
Also note that $R[x_1]/J$ is artinian by Lemma \ref{l3}(ii) and $R[x_1]/P_i$ is a simple artinian ring for every $i = 1,\dots,k$.
When $S/P_iS$ is considered as a ring of iterated differential polynomials, the appropriate versions of conditions (a), (b), (c) 
for $S/P_iS$ are satisfied and we may apply the induction. 
Therefore $L(n-1) = 0$ for every proper ideal $L$ of $S/P_iS$.

Repeating the argument from the first induction step, $P_i(1) = 0$ and hence also $P_iS(1) = 0$ for every $i = 1,\dots,k$. By Lemma \ref{l2}, either $I(n) = 0$ or
$I + P_iS = S$ for every $i = 1,\dots,k$. Assume the latter and let $P' := \cap_{i = 1}^k P_i$. Since $(P'+J)/J$ is the nilradical of $R[x_1]/J$,
there exists $m \in \mathbb{N}$ such that $P'^m \subseteq J$. Also note $P'S = \cap_{i = 1}^k P_iS$ and $(P'S)^m = P'^mS \subseteq I$.
Then $P'S + I /I$ is a nil ideal of $S/I$ and since we assume $P_iS + I = S$ for each $i = 1,\dots,k$, Lemma~\ref{l5} implies $I = S$.
This contradicts that $I$ is a proper ideal of $S$ and hence $P_iS + I \neq S$ for some $i$. Therefore $I(n) = 0$.  
\end{Proof}

\begin{Th} \label{t2}
Let $R$ be a right noetherian ring. Consider a ring of iterated differential polynomials 
$S = R[x_1,\delta_1,x_2,\delta_2,\dots,x_n,\delta_n]$ which satisfies 
\begin{enumerate}
\item[(a)] $\delta_i(R) \subseteq R$ for every $i = 1,\dots,n$,
\item[(b)] $1 \leq j < i \leq n$ there exist $u,v \in R[x_1,\dots,x_{j-1}]$ such that $\delta_i(x_j) = ux_j+v$, 
\item[(c)] if $J$ is an ideal of $R$ such that $\delta_i(J) \subseteq J$ for every $i = 1,\dots,n$, then either $J =0$ or $J = R$.
\end{enumerate}
Then $I(n) = 0$ for every proper ideal $I$ of $S$.
\end{Th}

\begin{Proof} 
Consider the ring $S_0 := R[x_1] \subseteq S$ and the set $\Delta := \{\delta_j|_{S_0} \mid 2 \leq j \leq n\}$ of derivations on $S_0$ induced by $\delta_2,\dots,\delta_n$.
Note that $0$ and $R$ are the only $\delta_1$-invariant ideals of $R$ which are also $\Delta$-invariant subsets of $S_0$ and $\partial(x_1)$ 
is a polynomial of degree at most one for every $\partial \in \Delta$. Hence Lemma \ref{l3} can be applied to $S_0$ and $\Delta$.

We proceed by induction on $n$ along the lines of the proof of Theorem \ref{t1}: If $n = 1$ and $I$ is a proper nonzero ideal of $S$, then $I$ is generated 
by a monic polynomial of $Z(S)$ of positive degree, hence $I(1) = 0$.

Now assume that the statement is true for rings with $n-1$ variables and consider a proper ideal $I \subseteq S = R[x_1,\delta_1,\dots,x_n,\delta_n]$.
Let $J: = I \cap S_0$. We distinguish two cases:

1) $J = 0$. Let $\Sigma$ be the set of all non-zero $\Delta$-invariant ideals of $S_0$. For each $M \in \Sigma$ let $f_M \in Z(S_0) \cap M$ 
be the monic polynomial such that $f_MS_0 = M$ (cf. Lemma \ref{l3}(i)). Set $X:=\{f_M \mid M \in \Sigma\}$.
Since a product of $\Delta$-invariant ideals of $S_0$ is $\Delta$-invariant, $X$ is a multiplicative set of regular elements of $S_0$ (observe $1 = f_{S_0}$) and 
since $X \subseteq Z(S_0)$, $X$ is a right denominator set of $S_0$. Then $X$ is a right denominator set of regular elements of $S$ and 
the canonical map  $\varphi\colon S \to S_0X^{-1}[x_2,\partial_2,\dots,x_n,\partial_n] = SX^{-1}$ given by 
$$\varphi(s) = s.1^{-1},s \in S_0,\ \varphi(x_i) = x_i, i = 2,\dots,n$$
is a right localization of $S$ with respect to $X$. Let us check, that the inductive assumption can be applied to 
the ring of iterated differential polynomials $S_0X^{-1}[x_2,\partial_2,\dots,x_n,\partial_n]$, i.e., appropriate versions of (a),(b) and (c) are satisfied
for $S_0X^{-1}[x_2,\partial_2,\dots,x_n,\partial_n]$.
Conditions (a) and (b) follow from the definition of derivations $\partial_2,\dots,\partial_n$. Let $L$ be a nonzero ideal of $S_0X^{-1}$ invariant under 
$\partial_2,\dots,\partial_n$. Let $L_0 = S_0 \cap L$. For $\ell \in L_0$ and $2\leq i \leq n$ is $(\delta_i(\ell)).1^{-1} = \partial_i(\ell.1^{-1}) \in L$. In particular, 
$\delta_i(\ell) \in L_0$, i. e., $L_0 \in \Sigma$. Then $f_{L_0} \in X$ and hence $L = L_0X^{-1} = S_0$.

 Since $SX^{-1}$ is right noetherian, \cite[Theorem~10.18(a)]{GW} implies $\varphi(I)X^{-1}$ is an ideal of $SX^{-1}$. The assumption $J = 0$ gives 
 $X \cap I = \emptyset$, i.e., $\varphi(I)X^{-1}$ is a proper ideal of $SX^{-1}$.
 The inductive assumption applied to  
 $SX^{-1}$ and Lemma \ref{l1} give that $\varphi(I(n-1)) = 0$. Observe that $X$ consists of regular elements of $S$, hence $\varphi$ is a monomorphism and 
 therefore $I(n-1) = 0$.

2) $J \neq 0$. By Lemma \ref{l3}(i) $J = fS_0$ for some monic polynomial from $Z(S_0)$. Consider the set $\Sigma$ of all $\Delta$-invariant ideals of $S_0$
which contain $J$ and let $\Sigma_m$ be the set of maximal elements from $\Sigma$. 
If $M \in \Sigma$, then $MS = M[x_2,x_3,\dots,x_n]$ is an ideal of $S$ 
and as noted in Section \ref{s2}, $S/MS$ has a structure of a ring of iterated differential polynomials,
$$S/MS \simeq (S_0/M)[x_2,\overline{\delta_2},\dots,x_n,\overline{\delta_n}]\,.$$
It follows directly from the construction in Section \ref{s2} that $(S_0/M)[x_2,\overline{\delta_2},\dots,x_n,\overline{\delta_n}]$
satisfies appropriate versions of (a) and (b). Since $M$ is a maximal element of $\Sigma$, is satisfies (c) as well.
By induction,  
$L(n-1)$ for every proper ideal $L$ of $S/MS$. By Lemma~\ref{l3}(iii), the set $\Sigma_m$ is finite and by Lemma~\ref{l3}(i), $M(1) = 0$ for every $M \in \Sigma_m$.
Therefore $MS(1) = 0$ for each $M \in \Sigma_m$.
Moreover, if $\Sigma_m = \{M_1,\dots,M_k\}$, then $J = M_1^{e_1}\cdots M_k^{e_k}$ for some $e_1,\dots,e_k \in \mathbb{N}$ and $(\cap \Sigma_m)^e \subseteq J$
for every $e \geq \sum_{i = 1}^k e_i$.
By Lemma~\ref{l2}, either $I(n) = 0$ or $I+MS = S$ for every $M \in \Sigma_m$. Assume the latter and let $N:= \cap_{i = 1}^k M_iS$. Then $N+I/I$ is a nil ideal 
of $S/I$ and, by Lemma~\ref{l5}, $I = S$. This contradiction proves that $I(n) = 0$.  
\end{Proof}

\begin{Rem} \label{r1}
Note that if $S$ is as in Theorem~\ref{t2}, then every nonzero ideal of $S$ contains a regular element. Indeed, keeping the notation from the proof of Theorem~\ref{t2},
we use the same inductive arguments: If $n = 1$, then $I$ contains a monic polynomial. 
If $n \geq 2$ and $J \neq 0$, then $J$ contains a monic polynomial. If $n \geq 2$ and $J = 0$, then, by induction, $\varphi(I)X^{-1}$ contains a regular element
$\varphi(f)\varphi(g)^{-1}$, where $f \in I$ and $g \in X$. But since $\varphi$ is a monomorphism and $\varphi(g)^{-1}$ is invertible in $SX^{-1}$, $f$ is a regular element of $S$.  
\end{Rem}

Let $L$ be a Lie algebra of finite dimension over a field $\mathbb{F}$. Recall that $L$ is {\em solvable} if it has a chain of subalgebras 
$0 = L_0 \subset L_1 \subset \cdots \subset L_n = L$ such that ${\rm dim}\ L_i = i$ and $L_i$ is an ideal in $L_{i+1}$ for each $i = 0,\dots,n-1$.
If there exists such a chain consisting of ideals in $L$, then we say that $L$ is {\em completely solvable}. If subalgebras in this chain satisfy $[L,L_i] \subseteq L_{i-1}$ for
$i = 1,2,\dots,n$, then $L$ is said to be {\em nilpotent}.

If $L$ is completely solvable, then its universal enveloping algebra $U(L)$ has a structure of a ring of iterated differential polynomials 
$U(L) = \mathbb{F}[x_1,\delta_1,\dots,x_n,\delta_n]$ (cf. \cite[p. xx]{GW}) which satisfies assumptions of Theorem~\ref{t2}. 
Combining Theorem~\ref{t2} and Proposition~\ref{p1}, we get 

\begin{Cor}
Let $L$ be a completely solvable Lie algebra of dimension $n$ over a field $\mathbb{F}$. Then $I(n) = 0$ for any proper ideal $I$ of $U(L)$. In particular, 
every projective $U(L)$-module is either finitely generated or free.  
\end{Cor}

Recall that, by \cite[Proposition~4.1]{PP}, $I(2) = 0$ if $I$ is a proper ideal of $U(L)$, where $L$ is a finite-dimensional solvable Lie algebra over a field 
of characteristic $0$. In the last part of this paper we show a similar result for ideals of finite codimension in $U(L)$, where $L$ is a finite-dimensional 
Lie algebra over an arbitrary field with $[L,L]$ nilpotent. If the characteristic of the field is zero, then $[L,L]$ is nilpotent for any solvable Lie algebra 
of finite dimension by \cite[Theorem~14.5.3]{MR}. The argument of \cite[Theorem~14.5.3]{MR} also shows $[L,L]$ is nilpotent if $L$ is completely solvable and the  
converse is true if the field is algebraically closed of characteristic zero.   

\begin{Lemma} \label{l4}
  Let $R$ be a right artinian ring, $\Delta$ a set of derivations on $R$ and let $\Sigma$ be the set of all $\Delta$-invariant ideals of $R$. Then the set $\Sigma_m$
	of all maximal elements of $\Sigma$ is finite and $\cap \Sigma_m$ is nilpotent.  
\end{Lemma}

\begin{Proof}
 Let $k$ be the length of $R_{R}$. If $M_1,\dots,M_{\ell}$ are pairwise different elements of $\Sigma_m$, then the Chinese reminder theorem implies that 
$\prod_{i = 1}^{\ell} R/M_i$ is a factor of $R$, hence $\ell \leq k$. Therefore $|\Sigma_m| \leq k$.

Let $J(R)$ be the Jacobson radical of $R$ and let $\pi\colon R \to R/J(R)$ be the canonical projection. We claim that for every maximal ideal $I$ of $R/J(R)$
there exists $M \in \Sigma_m$ such that $\pi(M) = I$. Indeed, let $I$ be a maximal ideal of a semisimple artinian ring $R/J(R)$. Then there exists an idempotent 
$e' \in R/J(R)$ such that $I = e' R/J(R)$. Since idempotents can be lifted modulo $J(R)$ (cf. \cite[Proposition~27.1]{AF}), there exists an idempotent $e \in R$ such that $e' = \pi(e)$.
Then $M' := ReR$ is $\Delta$-invariant and hence it is contained in some $M \in \Sigma_m$. Obviously $\pi(M')$ contains $I$, so either $\pi(M) = I$ or $\pi(M) = R/J(R)$.
The latter would mean $M + J(R) = R$ and consequently $M = R$ which is not the case. Therefore $\pi(M) = I$.

Finally let $J := \cap \Sigma_m$. Then, by the claim, $\pi(J)$ has to be contained in every maximal ideal of $R/J(R)$. In other words, $\pi(J) \subseteq J(R)$ which 
is a nilpotent ideal since $R$ is right artinian.    
\end{Proof}

Recall that an ideal $I$ of a right noetherian ring $R$ has the {\em right AR-property} if for every right ideal $A$ of $R$ there exists $m \in \mathbb{N}$
such that $A \cap I^m \subseteq AI$. A sequence $r_1,r_2,\dots,r_n \in R$ is said to be a {\em centralizing sequence} if for each $j \in \{1,\dots,n\}$ the image of $r_j$ in
$R/\sum_{i = 1}^{j-1} r_iR$ is a central element and $R \neq \sum_{i = 1}^n r_iR$. By \cite[Theorem~4.2.7(ii)]{MR}, any ideal of a right noetherian ring with 
a centralizing sequence of generators has the right AR-property.  

\begin{Prop} \label{p2}
Let $R$ be a right artinian ring. Consider a ring of iterated differential polynomials 
$S = R[x_1,\delta_1,x_2,\delta_2,\dots,x_n,\delta_n]$ which satisfies 
\begin{enumerate}
\item[(a)] $\delta_i(R) \subseteq R$ for every $i = 1,\dots,n$,
\item[(b)] $\delta_i(x_j) \in R$ for every $1 \leq j < i \leq n$,  
\item[(c)] if $J$ is an ideal of $R$ such that $\delta_i(J) \subseteq J$ for every $i = 1,\dots,n$, then either $J =0$ or $J = R$.
\end{enumerate}
If $I$ is a proper ideal of $S$ such that the right $R$-module $S/I$ has finite length, then $I$ is contained in a proper ideal of $S$ having the 
right AR-property. Moreover, $I(1) = 0$.
\end{Prop}

\begin{Proof}
Let $I$ be a proper ideal of $S$ such that $S/I$ is a right $R$-module of finite length. 
We claim that $I$ is contained in a proper ideal $K$ of $S$ which satisfies the right AR-property. 
More precisely, by induction on $n$ we show that  $I$ is contained in an ideal of $S$ which can be generated by a centralizing sequence.

If $n = 1$, then $I \neq 0$ and $I$ is generated by a monic polynomial $f \in Z(S)$ by Lemma~\ref{l3}(i). In this case we just set $K := I$.

Assume $n>1$ and let $S_0 := R[x_1]$, $J := I \cap S_0$. Note that $S_0/J$ is an $R$-submodule of $S/I$, hence $S_0/J$ is a right artinan ring. In particular, $J \neq 0$. 
Let $\Delta := \{\delta_i|_{S_0}\mid 2 \leq i \leq n \}$ be the set of derivations on $S_0$ induced by derivations $\delta_2,\dots,\delta_n$.
Further let $\Sigma$ be the set of all $\Delta$-invariant ideals of $S_0$ which contain $J$ and let $\Sigma_m$ be the set of all maximal elements of $\Sigma$.
By Lemma \ref{l3}(i), each $M \in \Sigma$ is generated by a unique monic polynomial $f_M \in Z(S_0)$ which is a non-zero polynomial of $M$ of the smallest possible degree. 
But regarding (b),  ${\rm deg}(\delta_i(f_M)) < {\rm deg} (f_M)$ for every $2 \leq i \leq n$. Hence $\delta_i(f_M) = 0$. Overall $f_M$ commutes not only with elements
of $S_0$ but also with $x_2,\dots,x_n$. In other words, $f_M \in Z(S)$.
   
By Lemma~\ref{l4} applied to $S_0/J$ and the set of derivations $\overline{\Delta}$ on $S_0/J$ which are induced by derivations of $\Delta$, the set $\Sigma_m$
is finite. Let $\Sigma_m = \{M_1,\dots,M_k\}$ and $L := \cap_{i = 1}^{k} M_i$. Then there exists $\ell \in \mathbb{N}$ such that $L^{\ell} \subseteq J$.   
For $i = 1,\dots,k$ let $M_i'$ be the ideal of $S$ generated by $M_i$. Then $L' := \cap_{i = 1}^k M_i' = LS$ and $L'^{\ell} = L^{\ell} S \subseteq I$.
In particular, $L'+I/I$ is a nil ideal of $S/I \neq 0$ and by Lemma~\ref{l5} there exists 
$i \in \{1,\dots,k\}$ such that $M_i'+I \neq S$.
 
Let $M \in \Sigma_m$ be such that $MS + I \neq S$. Recall that $M = f_M S_0$ for some $f_M$ in $Z(S)$. Using the notation introduced in Section \ref{s2},
let $$\pi \colon S \to (S_0/M)[x_2,\overline{\delta_2},\dots,x_n,\overline{\delta_n}] =: \overline{S} $$
be the homomorphism sending the term $rx_1^{e_1}x_2^{e_2}\cdots x_n^{e_n}$ to $(rx_1^{e_1}+M)x_2^{e_2} \cdots x_n^{e_n}$.
Note that ${\rm Ker}\ \pi = MS$ and $\pi(I)$ is a proper ideal of $\overline{S}$ and $\overline{S}/\pi(I)$ is finitely generated right 
$S_0/M$-module. Also the appropriate versions of (a), (b) and (c) are valid for 
$\overline{S}$: (a) and (b) are easy to check and (c) is satisfied since $M \in \Sigma_m$. By the induction, $\pi(I)$ is contained in $K'$ a proper 
ideal of $\overline{S}$ which is generated by a centralizing sequence $\pi(s_2),\dots,\pi(s_t)$ for some $s_2,\dots,s_t \in S$. Let $K:= \pi^{-1}(K')$, $s_1:=f_M$.
Then $s_1,s_2,\dots,s_{t}$ is a centralizing sequence of elements in $S$ generating $K$. Since $I \subseteq K$, we are done. 

By \cite[Theorem~2.2]{S}, $\cap_{i \in \mathbb{N}} K^i = \{s \in S \mid \exists a\in K $ such that $ 0 = s(1-a)\}$. If $K(1) = \cap_{i \in \mathbb{N}} K^{i} \neq 0$, then 
$K(1)$ contains a regular element $r$ by Remark \ref{r1}. But then $r(1-a) = 0$ for $a \in K$ implies $K = S$, which is not the case. Therefore $K(1) = 0$, hence also 
$I(1) = 0$.
\end{Proof}

\begin{Prop} \label{p3}
Let $L$ be a finite dimensional Lie algebra over any field such that $[L,L]$ is nilpotent. If $I$ is a proper ideal of finite codimension in the universal enveloping algebra 
$U(L)$, then $I(2) = 0$.
\end{Prop}

\begin{Proof}
Let $N:=[L,L]$ be the commutator of $L$ and let $R:=U(N)$ be the universal enveloping algebra of $N$. Let $x_1,\dots,x_k,x_{k+1},\dots,x_n$ be 
a basis of $L$ such that $x_1,\dots,x_k$ is a basis of $N$, then $S := U(L) = R[x_{k+1},\delta_{k+1},\dots,x_n,\delta_n]$ where $\delta_i(x_j) = [x_i,x_j] \in N 
\subseteq R$ for every $k \leq j < i \leq n$. Let $I$ be a proper ideal of $U(L)$ such that $U(L)/I$ is a finite dimensional algebra. 
Let $J:= R \cap I$. Observe that $J$ is an ideal of $R$ of finite codimension invariant under $\delta_{k+1},\dots,\delta_n$. In particular, $R/J$ is an 
artinian ring. Let $\Delta = \{\delta_i|_R \mid i = k+1 ,\dots, n\}$ be the set of derivations on $R$ given by restrictions of $\delta_{k+1},\dots,\delta_n$.
Consider the set $\Sigma$ of all $\Delta$-invariant ideals of $R$ which contain the ideal $J$. By Lemma~\ref{l4} applied to  $R/J$ and the set of 
derivations on $R/J$ induced by derivations from $\Delta$, we see that the set $\Sigma_m$ of maximal elements of $\Sigma$ is finite and there exists 
$m \in \mathbb{N}$ such that $(\cap \Sigma_m)^m \subseteq J$. Repeating the arguments from the proof of Proposition~\ref{p2}, there exists $M \in \Sigma_m$
such that $MS+I \neq S$: Let $\Sigma_m = \{ M_1,\dots,M_{\ell} \}$, $T:=\cap \Sigma_m$. Then $\cap_{i = 1}^{\ell} M_iS = TS$ and $(TS)^m \subseteq JS \subseteq I$, thus
$I + TS/I$ is a nil ideal of $S/I$. By Lemma~\ref{l5}, there exists $M \in \Sigma_m$ such that $MS + I \neq S$. 

Let $M \in \Sigma_m$ be such that $MS + I \neq S$. Using the notation of Section~\ref{s2}, the ring $\overline{S} := S/MS$ can be  seen as a ring 
of differential polynomials $\overline{S} = R/M[x_{k+1},\overline{\delta_{k+1}},\dots,x_{n},\overline{\delta_n}]$. Denote $\overline{I}:= I+MS/MS$ the 
image of $I$ in $S/MS$.
Let $\overline{\Delta} = \{\overline{\delta_i}|_{R/M} \mid i = k+1,\dots,n \}$ be the set of derivations on $R/M$ given by the restrictions of 
$\overline{\delta_{k+1}},\dots,\overline{\delta_n}$. Note that the only $\overline{\Delta}$-invariant ideals of $R/M$ are the trivial ones, therefore we
can apply Proposition~\ref{p2} to $\overline{I}$. That is, $\overline{I}(1) = 0$ and therefore $I(1) \subseteq MS$ by Lemma \ref{l1}.
   
Finally, recall that the Lie algebra $N$ is nilpotent and therefore, by \cite[Theorem~4.2]{M}, $R = U(N)$ is a right AR-ring, 
i.e., every right ideal of $R$ has the right AR-property. By \cite[Theorem~2.2]{S}, $M(1) = 0$ and also $MS(1)=0$. Then $I(2) \subseteq MS(1) = 0$.
\end{Proof}

\begin{Cor}
Let $V$ be a simple finite-dimensional module over a finite-dimensional Lie algebra with nilpotent commutator. Then ${\rm Ext}(V,V) \neq 0$.
\end{Cor}

\begin{Proof}
Consider $V$ as a (right) module over $U(L)$. Its annihilator $I$ is an ideal of $U(L)$ of finite codimension. If ${\rm Ext}(V,V) = 0$, then 
arguments from the proof of \cite[Corollary~3.8]{PP} show that $I = I^2$ which contradicts Proposition~\ref{p3}.
\end{Proof}
\bibliographystyle{abbrv}
\bibliography{references}

\begin{thebibliography}{1}

\bibitem{AF}
F.~W. Anderson and K.~R. Fuller.
\newblock {\em Rings and categories of modules}, volume Vol. 13 of {\em
  Graduate Texts in Mathematics}.
\newblock Springer-Verlag, New York-Heidelberg, 1974.

\bibitem{B}
H.~Bass.
\newblock Big projective modules are free.
\newblock {\em Illinois J. Math.}, 7:24--31, 1963.

\bibitem{GW}
K.~R. Goodearl and R.~B. Warfield, Jr.
\newblock {\em An introduction to noncommutative {N}oetherian rings}, volume~61
  of {\em London Mathematical Society Student Texts}.
\newblock Cambridge University Press, Cambridge, second edition, 2004.

\bibitem{M}
J.~C. McConnell.
\newblock The intersection theorem for a class of non-commutative rings.
\newblock {\em Proc. London Math. Soc. (3)}, 17:487--498, 1967.

\bibitem{MR}
J.~C. McConnell and J.~C. Robson.
\newblock {\em Noncommutative {N}oetherian rings}, volume~30 of {\em Graduate
  Studies in Mathematics}.
\newblock American Mathematical Society, Providence, RI, revised edition, 2001.
\newblock With the cooperation of L. W. Small.

\bibitem{PP}
P.~P\v{r}\'{\i}hoda and G.~Puninski.
\newblock Iterated power intersections of ideals in rings of iterated
  differential polynomials.
\newblock {\em J. Algebra Appl.}, 12(7):1350020, 10, 2013.

\bibitem{S}
P.~F. Smith.
\newblock The {A}rtin-{R}ees property.
\newblock In {\em Paul {D}ubreil and {M}arie-{P}aule {M}alliavin {A}lgebra
  {S}eminar, 34th {Y}ear ({P}aris, 1981)}, volume 924 of {\em Lecture Notes in
  Math}, pages pp 197--240. Springer, Berlin-New York, 1982.

\end{thebibliography}
\end{document}